\documentclass[11pt]{amsart}

\usepackage{amsmath,amssymb}
\usepackage{mathrsfs}
\usepackage{xcolor}
\usepackage{hyperref}
\newtheorem{theorem}{Theorem}
\newtheorem{proposition}{Proposition}
\newtheorem{lemma}{Lemma}
\theoremstyle{definition}
\newtheorem{question}{Question}
\theoremstyle{remark}
\newtheorem{remark}{Remark}

\DeclareMathOperator{\GL}{GL}
\DeclareMathOperator{\Lip}{Lip}
\DeclareMathOperator{\dist}{dist}
\DeclareMathOperator{\spanop}{span}
\newcommand{\F}{\mathbb{F}}
\newcommand{\R}{\mathbb{R}}
\newcommand{\C}{\mathbb{C}}
\newcommand{\N}{\mathbb{N}}
\newcommand{\Nzero}{\mathbb{N}_{0}}
\newcommand{\Tcircle}{\mathbb{T}}
\newcommand{\K}{\mathcal{K}}
\newcommand{\Lcal}{\mathcal{L}}
\newcommand{\AHX}{\mathfrak{X}_{K}}

\title[Strong Structural Stability]
{Strong Structural Stability of Scalar-\allowbreak Plus-\allowbreak Compact
Operators on Banach Spaces}
\author{Xinxing Wu}
\address{School of Mathematics and Statistics, Guizhou University of Finance and Economics, Guiyang, Guizhou 550025, China}
\email{wuxinxing5201314@163.com}
\author{Mengjuan Wang}
\address{School of Mathematics and Statistics, Guizhou University of Finance and Economics, Guiyang, Guizhou 550025, China}
\email{wangmengjuan@mail.gufe.edu.cn}
\author{Xu Zhang}
\address{Department of Mathematics\\
Shandong University, Weihai, Shandong, 264209,  China}
\email{xu\_zhang\_sdu@mail.sdu.edu.cn (Corresponding author)}

\begin{document}

\begin{abstract}
We prove that every strongly structurally stable operator of the form
\(\lambda I+K\) on a real or complex Banach space, where \(\lambda\) is a scalar
and \(K\) is compact, is
hyperbolic. Applying this result to the real Banach space constructed by
Argyros and Haydon and its complexification, both of which are
infinite-dimensional separable Banach spaces with the scalar-plus-compact
property, we show that neither space supports a nonhyperbolic and
strongly structurally stable operator, thereby giving a negative answer to the
part of the open problem of
Bernardes and Messaoudi \cite{BM20} concerning strong structural stability.
\end{abstract}

\maketitle

\section{Introduction}

Structural stability was introduced by Andronov and Pontryagin \cite{AP37} and
subsequently became a central concept in dynamical systems
\cite{GuckenheimerHolmes83,KatokHasselblatt95,Shub87}. In the differentiable setting the
perturbations are usually measured in a \(C^{1}\) norm, whereas the
operator-theoretic formulation used by Pugh \cite{Pugh69} and by Bernardes and
Messaoudi \cite{BM20} allows bounded Lipschitz perturbations with small
supremum norm and small Lipschitz constant. A closely related result is the classical
Grobman--Hartman theorem \cite{Grobman59,Hartman60,Hartman63}.
In its discrete-time form, it asserts that the
dynamics near a hyperbolic fixed point of a finite-dimensional \(C^{1}\) diffeomorphism are
locally topologically conjugate to the dynamics of its derivative.
Its extension to Banach spaces is due
independently to Palis \cite{Palis68} and Pugh \cite{Pugh69}.

For invertible linear operators on Banach spaces,
Pugh~\cite{Pugh69} established that hyperbolicity implies structural
stability, using the characterization of hyperbolicity by a direct sum
decomposition into two closed invariant subspaces on which the operator and its
inverse, respectively, are strict contractions after an equivalent renorming.
Conversely, Robbin \cite{Robbin72} proved that every
structurally stable invertible linear operator on a finite-dimensional space
is hyperbolic. In infinite dimensions, however, structural stability need not imply
hyperbolicity. Indeed, Bernardes and Messaoudi \cite{BM20,BM21}
constructed invertible bilateral weighted backward shifts on
\(\ell_{p}(\mathbb{Z})\), \(1\leq p<\infty\), and on \(c_{0}(\mathbb{Z})\) that are
structurally stable and not hyperbolic.

Bernardes et al.~\cite{BCDMP18} introduced a notion now known as
generalized hyperbolicity and proved that every generalized hyperbolic operator
has the shadowing property.
Generalized hyperbolicity is a weaker condition than hyperbolicity:
this class contains every invertible hyperbolic operator, as well as the
nonhyperbolic weighted backward shifts described above \cite{BM20}.
Cirilo et al.~\cite{CGP21} later showed that the class of generalized hyperbolic
operators is open in the operator norm topology.
Regarding structural stability, Bernardes and Messaoudi~\cite{BM20}
established that every generalized hyperbolic operator on a
Banach space is strongly structurally stable, with the conjugacy uniformly
close to the identity, and consequently extended the Grobman--Hartman theorem
to generalized hyperbolic fixed points of \(C^{1}\) diffeomorphisms on Banach
spaces. They also
posed as an open problem whether strong structural
stability is equivalent to ordinary structural stability.

In light of the existence of nonhyperbolic and strongly structurally
stable weighted shifts, Bernardes and Messaoudi~\cite{BM20} posed the following
question:

\begin{question}[{\cite[p.~4359]{BM20}}]\label{q:bm}
Does every infinite-dimensional (separable) Banach space support a nonhyperbolic
\((\)strongly\()\) structurally stable operator?
\end{question}

Throughout this paper, let \(\F\in\{\R,\C\}\),
\(\N=\{1,2,\ldots\}\), \(\Nzero=\{0,1,2,\ldots\}\), and
\(\Tcircle=\{z\in\C:|z|=1\}\). For a Banach
space \(X\) over \(\F\), let \(\Lcal(X)\) denote the set of all bounded
\(\F\)-linear operators on \(X\), \(\K(X)\) the set of compact
operators on \(X\), \(\GL(X)\) the set of invertible elements of
\(\Lcal(X)\), and \(I_{X}\) the identity operator on \(X\).
For a real Banach space \(X\), let
\(X_{\C}=X\oplus\mathbf{i}X\) be its complexification.
For \(T\in\Lcal(X)\), define its complexification
\(T_{\C}:X_{\C}\to X_{\C}\) by
\[
T_{\C}(x+\mathbf{i}y)=T(x)+\mathbf{i}T(y)\quad(x,y\in X).
\]
We equip \(X_{\C}\) with the norm
\[
\|x+\mathbf{i}y\|_{X_{\C}}
:=\sup_{\theta\in[0,2\pi]}\|x\cos(\theta)-y\sin(\theta)\|_{X}
\quad(x,y\in X).
\]
This norm extends the norm on \(X\) and is invariant under complex
conjugation~\cite[Section~2]{MST99}.

According to Argyros and Haydon~\cite{AH11}, a Banach space
\(X\) over \(\F\) has the scalar-plus-compact property if
\(\Lcal(X)=\F I_{X}+\K(X),\)
where
\(\F I_{X}+\K(X)=\{\lambda I_{X}+K:\lambda\in\F,\ K\in\K(X)\}.\)

The following theorem is the main result. It shows that strong structural
stability implies hyperbolicity for compact perturbations of scalar multiples
of the identity.

\begin{theorem}\label{thm:main}
Let \(X\) be a Banach space over \(\F\), and let
\(T\in\GL(X)\) satisfy \(T=\lambda I_{X}+K\),
where \(\lambda\in\F\) and \(K\in\K(X)\).
If \(T\) is strongly structurally stable, then it is hyperbolic.
\end{theorem}

In particular, Theorem~\ref{thm:main} implies that every
strongly structurally stable invertible operator on a
Banach space with the scalar-plus-compact property is hyperbolic.
For example, Argyros and Haydon~\cite[Theorem~7.4]{AH11} constructed an
infinite-dimensional separable real Banach space \(\AHX\) with this
property. Moreover, Theorem~\ref{thm:ah-complexification}
shows that the complexification \((\AHX)_{\C}\) of \(\AHX\) is
an infinite-dimensional separable complex Banach space with the
scalar-plus-compact property. Theorem~\ref{thm:main} therefore implies that
neither \(\AHX\) nor \((\AHX)_{\C}\) supports a nonhyperbolic and strongly
structurally stable operator. Thus the part of Question~\ref{q:bm}
concerning strong structural stability has a negative answer
over both real and complex scalar fields.

\section{Preliminaries}

For \(r>0\), set
\(B_{X}(0,r)=\{x\in X:\|x\|\leq r\}.\)
The spectrum of a bounded \(\C\)-linear operator \(T\) is denoted by \(\sigma(T)\).
For a bounded \(\R\)-linear operator \(T\), we set
\(\sigma(T):=\sigma(T_{\C})\).
For a nonempty subset \(M\subseteq X\), set
\[
\dist(x,M)=\inf\{\|x-y\|:y\in M\}\quad(x\in X).
\]
For \(E\subseteq X\), \(\spanop E\) denotes the linear span of \(E\).

We recall the following deep result of Argyros and Haydon~\cite{AH11}.

\begin{theorem}[{\cite[Theorem 7.4]{AH11}}]
\label{thm:argyros-haydon}
There exists an infinite-dimensional separable real Banach space \(\AHX\) such
that, for every \(T\in\Lcal(\AHX)\), there exists a scalar \(\lambda\in\R\)
such that \(T-\lambda I_{\AHX}\in\K(\AHX).\)
\end{theorem}

Theorem~\ref{thm:argyros-haydon} gives \(\Lcal(\AHX)\subseteq
\K(\AHX)+\R I_{\AHX}\). The reverse inclusion follows directly
from \(I_{\AHX}\in\Lcal(\AHX)\) and \(\K(\AHX)\subseteq\Lcal(\AHX)\).
Consequently,
\[
\Lcal(\AHX)=\K(\AHX)\oplus\R I_{\AHX}.
\]

B{\"u}hler and Salamon~\cite[Example~5.2.9]{BS18} state that the
complexification \((\AHX)_{\C}\) of the real Banach space \(\AHX\) also
has the scalar-plus-compact property.

\begin{theorem}[{\cite[Example~5.2.9]{BS18}}]
\label{thm:ah-complexification}
For the real Banach space \(\AHX\) of Theorem \ref{thm:argyros-haydon}, the
complexification \((\AHX)_{\C}\) is an infinite-dimensional separable complex
Banach space satisfying
\(\Lcal((\AHX)_{\C})=\K((\AHX)_{\C})\oplus\C I_{(\AHX)_{\C}}.\)
\end{theorem}

For a map \(f:E\to Y\), where \(E\) is a nonempty subset of a
Banach space \(X\) over \(\F\) and \(Y\) is a normed space, define
\[
\|f\|_{\infty}=\sup_{x\in E}\|f(x)\|
\quad\text{and}\quad
\Lip(f)=\sup_{x,y\in E,\, x\neq y}\frac{\|f(x)-f(y)\|}{\|x-y\|}.
\]
Both quantities are allowed to be \(+\infty\); if \(E\) is a
singleton, we set \(\Lip(f)=0\).

Let \(X\) be a complex Banach space and \(T\in\Lcal(X)\). A point
\(\mu\in\sigma(T)\) is an isolated spectral value of \(T\) if \(\mu\)
is an isolated point of \(\sigma(T)\). For such \(\mu\), the Riesz
projection associated with \(\mu\) is
\[
P_{\mu}=\frac{1}{2\pi\mathbf{i}}\int_{\Gamma} (zI_{X}-T)^{-1}\,dz,
\]
where \(\Gamma\) is a positively oriented rectifiable simple closed curve in
\(\C\setminus\sigma(T)\) enclosing \(\mu\) but no other point of \(\sigma(T)\)
\cite[Chapter~III, \S\S~6.4--6.5]{Kato95}.
Thus \(P_{\mu}\) is a projection on \(E_{\mu}=P_{\mu}(X)\) along
\(F_{\mu}=(I_{X}-P_{\mu})(X)\)
\cite[Chapter~III, Theorem~6.17]{Kato95}.
The subspace \(E_{\mu}\) is called the spectral subspace of \(T\) associated
with \(\mu\) \cite[Section~V.1]{GGK90}.
By \cite[Chapter~III, Theorem~6.17]{Kato95}, \(P_{\mu}\) commutes with
\(T\) (i.e., \(P_{\mu}\circ T=T\circ P_{\mu}\)), \(E_{\mu}\) is \(T\)-invariant, and
\(\sigma(T|_{E_{\mu}})=\{\mu\}\).
If \(E_{\mu}\) is finite-dimensional, then \(\mu\) is an eigenvalue of \(T\).
In this case \(\dim E_{\mu}\) is called the (algebraic) multiplicity of the
eigenvalue \(\mu\) of \(T\)~\cite[Chapter~III, \S~6.5]{Kato95}.

An operator \(T\in\GL(X)\) is hyperbolic when
\(
\sigma(T)\cap\Tcircle=\varnothing .
\)
According to Bernardes and Messaoudi~\cite{BM20}, the operator
\(T\) is strongly structurally stable if, for every
\(\gamma>0\), there exists \(\varepsilon>0\) such that,
for every bounded Lipschitz map \(\beta:X\to X\) satisfying
\[
\|\beta\|_{\infty}\leq\varepsilon \quad\text{and}\quad \Lip(\beta)\leq\varepsilon ,
\]
there exists a homeomorphism \(H:X\to X\) satisfying
\[
H\circ T=(T+\beta)\circ H \quad\text{and}
\quad \|H-I_{X}\|_{\infty}\leq\gamma .
\]
In this case, \(T\) and \(T+\beta\) are said to be topologically conjugate through \(H\).

The bound on \(\|H-I_{X}\|_{\infty}\) in the definition of strong structural
stability yields the following comparison between the forward orbits of
\(T+\beta\) and those of \(T\).

\begin{lemma}\label{lem:orbit-condition}
Let \(X\) be a Banach space over \(\F\), and let
\(T\in\GL(X)\) be strongly structurally stable. Then, for every
\(\gamma>0\), there exists \(\varepsilon>0\) such that,
for every bounded Lipschitz map \(\beta:X\to X\) with
\begin{equation}\label{two-estimates}
\|\beta\|_{\infty}\leq\varepsilon \quad\text{and}\quad \Lip(\beta)\leq\varepsilon ,
\end{equation}
and for every \(y\in X\), there exists \(x\in X\) satisfying
\[
\|(T+\beta)^{n}(y)-T^{n}(x)\|\leq\gamma\quad(n\in\Nzero).
\]
\end{lemma}

\begin{proof}
Fix \(\gamma>0\), and choose \(\varepsilon>0\) given by
strong structural stability for this \(\gamma\). For every bounded
Lipschitz map \(\beta:X\to X\) satisfying condition~\eqref{two-estimates},
strong structural stability gives a homeomorphism \(H:X\to X\) such that
\(H\circ T=(T+\beta)\circ H\) and \(\|H-I_{X}\|_{\infty}\leq\gamma.\)
For every \(y\in X\), set \(x=H^{-1}(y)\).
Iterating \(H\circ T=(T+\beta)\circ H\) gives
\[
(T+\beta)^{n}(y)=(T+\beta)^{n}(H(x))=H(T^{n}(x))\quad(n\in\Nzero).
\]
Using \(\|H-I_{X}\|_{\infty}\leq\gamma\), one obtains
\[
\|(T+\beta)^{n}(y)-T^{n}(x)\|
=\|(H-I_{X})(T^{n}(x))\|
\leq\|H-I_{X}\|_{\infty}\leq\gamma\
(n\in\Nzero).
\]
\end{proof}

The following standard consequence of the contraction mapping theorem
ensures that the maps \(T+\beta\) constructed below are homeomorphisms.
Nevertheless, the proof is included for completeness.

\begin{lemma}\label{lem:homeomorphism}
Let \(X\) be a Banach space over \(\F\), let \(T\in\GL(X)\) and \(\beta:X\to X\) be Lipschitz. If
\(
\|T^{-1}\|\cdot \Lip(\beta)<1,
\)
then \(T+\beta\) is a homeomorphism of \(X\).
\end{lemma}

\begin{proof}
Since \(T\in\GL(X)\), \(T\) is a homeomorphism. The identity
\(T+\beta=T\circ(I_{X}+T^{-1}\circ\beta)\)
therefore shows that it suffices to prove that \(I_{X}+T^{-1}\circ\beta\) is a homeomorphism.
For \(y\in X\), the equation
\[
(I_{X}+T^{-1}\circ\beta)(x)=x+T^{-1}(\beta(x))=y
\]
is equivalent to the fixed point equation
\[
x=y-T^{-1}(\beta(x)).
\]

Given any fixed \(y\in X\), the map \(x\mapsto y-T^{-1}(\beta(x))\) is a contraction with
constant at most \(\|T^{-1}\|\cdot \Lip(\beta)<1\).
The contraction mapping theorem gives a unique solution
of \(x=y-T^{-1}(\beta(x))\). Hence \(I_{X}+T^{-1}\circ\beta\) is bijective.
Moreover, \(I_{X}+T^{-1}\circ\beta\) is Lipschitz, hence continuous, because
\(T^{-1}\circ\beta\) is Lipschitz.

It remains to prove that the inverse of \(I_{X}+T^{-1}\circ\beta\) is continuous. Let \(y_{1},y_{2}\in X\) be
arbitrary, and set
\[
x_{1}=(I_{X}+T^{-1}\circ\beta)^{-1}(y_{1})
\quad\text{and}\quad
x_{2}=(I_{X}+T^{-1}\circ\beta)^{-1}(y_{2}).
\]
Then
\[
x_{1}+T^{-1}(\beta(x_{1}))=y_{1}
\quad\text{and}\quad
x_{2}+T^{-1}(\beta(x_{2}))=y_{2}.
\]
Subtracting the two identities gives
\[
x_{1}-x_{2}+T^{-1}(\beta(x_{1})-\beta(x_{2}))=y_{1}-y_{2},
\]
implying
\[
\|x_{1}-x_{2}\|
\leq \|y_{1}-y_{2}\|+\|T^{-1}\|\cdot
\Lip(\beta)\cdot \|x_{1}-x_{2}\|.
\]
This, together with \(1-\|T^{-1}\|\cdot\Lip(\beta)>0\), gives
\[
\|(I_{X}+T^{-1}\circ\beta)^{-1}(y_{1})
-(I_{X}+T^{-1}\circ\beta)^{-1}(y_{2})\|=
\|x_{1}-x_{2}\|\leq\frac{\|y_{1}-y_{2}\|}{1-\|T^{-1}\|\cdot \Lip(\beta)}.
\]
Therefore \((I_{X}+T^{-1}\circ\beta)^{-1}\) is Lipschitz, and hence continuous.
\end{proof}

\begin{lemma}\label{lem:compact-direction}
Let \(X\) be an infinite-dimensional Banach space over \(\F\), \(F\subseteq X\) a
finite-dimensional subspace of \(X\), and \(K\in\K(X)\). For every
\(\alpha>0\), there exists \(e\in X\) with \(\|e\|=1\)
such that
\(\dist(e,F)\geq\frac{1}{4}\) and \(\|K(e)\|<\alpha\).
\end{lemma}

\begin{proof}
Fix \(\alpha >0\).
Riesz's lemma gives unit vectors \(x_{n}\) such that
\begin{equation}\label{eq:compact-direction-distance}
\dist(x_{n},F+\spanop\{x_{1},\ldots,x_{n-1}\})>\frac{3}{4}
\quad(n\in\N).
\end{equation}
Since \(K\) is compact, the sequence \((K(x_{n}))_{n\in\N}\) has a convergent
subsequence. Hence there exist \(m,n\in\N\) with \(m<n\) such that
\begin{equation}\label{eq:compact-direction-K}
\|K(x_{n})-K(x_{m})\|<\frac{\alpha}{2}.
\end{equation}
Since \(m<n\), we have
\(x_{m}\in \operatorname{span}\{x_{1},\ldots,x_{n-1}\}
\subseteq F+\operatorname{span}\{x_{1},\ldots,x_{n-1}\},\)
so inequality~\eqref{eq:compact-direction-distance} gives
\(\|x_{n}-x_{m}\|>\frac{3}{4}.\)
Thus the vector
\(e=\frac{x_{n}-x_{m}}{\|x_{n}-x_{m}\|}\)
is well defined and satisfies \(\|e\|=1\). Using \(\|x_{n}-x_{m}\|>\frac{3}{4}\) and
\eqref{eq:compact-direction-K}, we obtain
\[
\|K(e)\|
=\frac{\|K(x_{n})-K(x_{m})\|}{\|x_{n}-x_{m}\|}
<\frac{\alpha/2}{3/4}
<\alpha .
\]

To estimate \(\dist(e,F)\), fix any \(y\in F\). Since
\(\|x_{n}-x_{m}\|\cdot y+x_{m}\in F+\operatorname{span}\{x_{1},\ldots,x_{n-1}\},\)
inequality \eqref{eq:compact-direction-distance} and
\(\|x_{n}-x_{m}\|\leq \|x_{n}\|+\|x_{m}\|=2\) give
\[
\|e-y\|
=\frac{\|x_{n}-(x_{m}+\|x_{n}-x_{m}\|\cdot y)\|}{\|x_{n}-x_{m}\|}
>\frac{3/4}{2}
=\frac{3}{8},
\]
implying
\(\dist(e,F)\geq\frac{3}{8}\geq\frac{1}{4}\).
\end{proof}

\begin{lemma}\label{lem:interpolation}
Let \(X\) be a Banach space over \(\F\). Fix \(N\in\N\),
\(R\geq0\), \(\mathbb{K}\in\{\R,\C\}\), and
\(a_{0},\ldots,a_{N}\in X\). Suppose that
\(g:\{a_{0},\ldots,a_{N}\}\to\mathbb{K}\)
is Lipschitz and satisfies
\(|g(a_{m})|\leq R\) \((m\in\{0,\ldots,N\}).\)
Then there exists a bounded Lipschitz map \(G:X\to\mathbb{K}\) such that
\[
G(a_{m})=g(a_{m})\quad(m\in\{0,\ldots,N\}),
\]
and
\[
\|G\|_{\infty}\leq R,\quad \Lip(G)\leq2\Lip(g).
\]
\end{lemma}

\begin{proof}
Since \(\mathbb{K}\subseteq\C\), applying the McShane
extension theorem to the real and imaginary parts of \(g\) gives extensions
\(u,v:X\to\R\) without increasing their respective Lipschitz constants.
We take \(v\equiv0\) when \(\mathbb{K}=\R\).
In particular, \(\Lip(u),\Lip(v)\leq\Lip(g)\) and
\[
u(a_{m})=\operatorname{Re}(g(a_{m})),\quad
v(a_{m})=\operatorname{Im}(g(a_{m}))\quad(m\in\{0,\ldots,N\}).
\]
Set \(\tilde{G}=u+\mathbf{i}v:X\to\mathbb{K}\). Then
\(\tilde{G}\) extends \(g\), and for any \(x,y\in X\), we have
\[
|\tilde{G}(x)-\tilde{G}(y)|
\leq |u(x)-u(y)|+|v(x)-v(y)|
\leq2\cdot \Lip(g)\cdot\|x-y\|.
\]
Thus \(\Lip(\tilde{G})\leq2\Lip(g)\).

To obtain a bounded extension, define
\(P:\mathbb{K}\to\{z\in\mathbb{K}:|z|\leq R\}\) by
\[
P(z)=
\begin{cases}
z,& |z|\leq R,\\
\frac{Rz}{|z|},& |z|>R.
\end{cases}
\]
This map is the metric projection onto
\(\{z\in\mathbb{K}:|z|\leq R\}\) in the real Hilbert space \(\mathbb{K}\).
It is therefore nonexpansive, that is,
\[
|P(x_{1})-P(x_{2})|\leq|x_{1}-x_{2}|\quad(x_{1}, x_{2}\in\mathbb{K}).
\]
Since \(\tilde{G}(a_{m})=g(a_{m})\) and \(|g(a_{m})|\leq R\), we have
\(P(\tilde{G}(a_{m}))=\tilde{G}(a_{m})\) for \(m\in\{0,\ldots,N\}\).
Thus the map \(G=P\circ\tilde{G}\) extends \(g\), satisfies
\(\|G\|_{\infty}\leq R\), and satisfies
\(\Lip(G)\leq\Lip(\tilde{G})\leq2\Lip(g)\)
 because \(P\) is nonexpansive.
\end{proof}

\begin{lemma}\label{lem:radial}
Fix \(\nu\in\Tcircle\), \(\delta>0\), and \(N\in\N\).
For every \(R>0\), the assignment
\[
\nu^{m}(R+m\delta)\mapsto \nu^{m+1}\quad(m\in\{0,\ldots,N\})
\]
is well defined and Lipschitz with constant at most \(\frac{1}{R}\).
\end{lemma}

\begin{proof}
Fix \(R>0\). Since \(|\nu^{m}(R+m\delta)|=R+m\delta\) and \(\delta>0\), the points
\(\nu^{m}(R+m\delta)\), \(m\in\{0,\ldots,N\}\), are pairwise distinct. Hence the
assignment is well defined.
For any \(m,n\in\{0,\ldots,N\}\) with \(m<n\), using
\(|\nu|=1\), we have
\[
\begin{aligned}
&|\nu^{n}(R+n\delta)-\nu^{m}(R+m\delta)|^{2}\\
&\qquad=[\nu^{n}(R+n\delta)-\nu^{m}(R+m\delta)]
\cdot\overline{[\nu^{n}(R+n\delta)-\nu^{m}(R+m\delta)]}\\
&\qquad=[\nu^{n}(R+n\delta)-\nu^{m}(R+m\delta)]
\cdot [\nu^{-n}(R+n\delta)-\nu^{-m}(R+m\delta)]\\
&\qquad=(R+n\delta)^{2}+(R+m\delta)^{2}
-(R+n\delta)\cdot(R+m\delta)\cdot(\nu^{n-m}+\nu^{m-n})\\
&\qquad=(n-m)^{2}\cdot \delta^{2}
  +(R+n\delta)\cdot(R+m\delta)\cdot(\nu^{n}-\nu^{m})\cdot (\nu^{-n}-\nu^{-m})\\
&\qquad=(n-m)^{2}\cdot \delta^{2}
  +(R+n\delta)\cdot(R+m\delta)\cdot|\nu^{n}-\nu^{m}|^{2}\\
&\qquad\geq R^{2}\cdot |\nu^{n}-\nu^{m}|^{2}.
\end{aligned}
\]
This, together with
\(|\nu^{n+1}-\nu^{m+1}|=|\nu^{n}-\nu^{m}|\), implies
\[
\frac{|\nu^{n+1}-\nu^{m+1}|}
{|\nu^{n}(R+n\delta)-\nu^{m}(R+m\delta)|}
=\sqrt{\frac{|\nu^{n+1}-\nu^{m+1}|^{2}}
{|\nu^{n}(R+n\delta)-\nu^{m}(R+m\delta)|^{2}}}
\leq\frac{1}{R}.
\]
\end{proof}

\section{Main results}

We prove Theorem \ref{thm:main} using two propositions.
Proposition \ref{prop:essential-scalar} treats compact perturbations of scalar
multiples of the identity on infinite-dimensional spaces, where the scalar
has modulus one.
Proposition \ref{prop:finite-spectrum} treats operators \(T\) with an isolated
spectral value on the unit circle and a finite-dimensional associated
spectral subspace (of \(T_{\C}\) in the real case).

\begin{proposition}\label{prop:essential-scalar}
Let \(X\) be an infinite-dimensional Banach space over \(\F\),
and let \(T\in\GL(X)\) satisfy
\(T=\lambda I_{X}+K,\)
where \(\lambda\in\F\), \(|\lambda|=1\), and \(K\in\K(X)\).
Then \(T\) is not strongly structurally stable.
\end{proposition}

\begin{proof}
Suppose, on the contrary, that \(T\) is strongly structurally stable. Let
\(\varepsilon_{0}>0\) be given by Lemma \ref{lem:orbit-condition} with
\(\gamma=1\), and set
\[
\varepsilon:=\frac{1}{2}\min\left\{1,\|T^{-1}\|^{-1},\varepsilon_{0}\right\},\quad
\delta:=\frac{\varepsilon}{32},
\quad\text{and}\quad N:=\left\lceil\frac{16}{\delta}\right\rceil.
\]
Then \(N\delta\geq16\).
Since \(\overline{K(B_{X}(0,N\delta+2))}\) is compact, there exist
\(p_{1},\ldots,p_{s}\in X\) such that
\(
\overline{K(B_{X}(0,N\delta+2))}
\subseteq\bigcup_{j=1}^{s}\{y\in X:\|y-p_{j}\|<\tfrac{\delta}{64}\}.
\)
Let
\(F=\spanop\{p_{1},\ldots,p_{s}\}.\)
Then \(F\) is finite-dimensional and
\begin{equation}\label{eq:distance-KX-F}
\dist(K(x),F)\leq\frac{\delta}{64}\quad(\|x\|\leq N\delta+2).
\end{equation}
By Lemma~\ref{lem:compact-direction}, there exists \(e\in X\)
with \(\|e\|=1\) such that
\[
\dist(e,F)\geq\frac{1}{4}\quad \text{and} \quad
\|K(e)\|<
\frac{\delta}{1+N\delta}.
\]

Let
\[
z_{m}=\lambda^{m}(1+m\delta)
\quad\text{and}\quad
a_{m}=z_{m} e \quad(m\in\{0,\ldots,N\}).
\]
By Lemma~\ref{lem:radial}, applied with \(\nu=\lambda\) and \(R=1\), the assignment
\[
z_{m}\mapsto\lambda^{m+1}\quad(m\in\{0,\ldots,N\})
\]
has Lipschitz constant at most \(1\).
This, together with \(\|a_{m}-a_{n}\|=|z_{m}-z_{n}|\),
implies that the assignment \(a_{m}\mapsto\lambda^{m+1}\) is well
defined and has Lipschitz constant at most \(1\).
Since \(|\lambda|=1\), for every \(m\in\{0,\ldots,N\}\), we have
\(|\lambda^{m+1}|=1\) and \(|z_{m}|=1+m\delta\leq1+N\delta.\)
Lemma~\ref{lem:interpolation}, applied with \(R=1\) and \(R=1+N\delta\),
respectively, therefore gives bounded Lipschitz maps
\(\varphi,\psi:X\to\F\) such that
\[
\varphi(a_{m})=\lambda^{m+1},\quad \psi(a_{m})=z_{m}
\quad(m\in\{0,\ldots,N\}),
\]
\[
\|\varphi\|_{\infty}\leq1,\quad \Lip(\varphi)\leq2,
\]
and
\[
\|\psi\|_{\infty}\leq 1+N\delta,\quad \Lip(\psi)\leq2.
\]

Define
\(\beta(x)=\delta\varphi(x)e-\psi(x)K(e) .\)
The interpolation identities for \(\varphi\) and \(\psi\), together with
\(z_{m+1}-\lambda z_{m}=\delta\lambda^{m+1}\), \(a_{m}=z_{m}e\), and
\(T=\lambda I_{X}+K\), give
\[
\begin{aligned}
\beta(a_{m})
&=\delta\lambda^{m+1}e-z_{m}K(e)\\
&=(z_{m+1}-\lambda z_{m})e-z_{m}K(e)\\
&=a_{m+1}-T(a_{m})\quad(m\in\{0,\ldots,N-1\}).
\end{aligned}
\]
Thus the map \(S=T+\beta\) satisfies
\[
S(a_{m})=a_{m+1}\quad(m\in\{0,\ldots,N-1\}).
\]
Moreover, the bounds on
\(\varphi,\psi\), and \(K(e)\) give
\[
\|\beta\|_{\infty}
\leq \delta+(1+N\delta)\cdot \|K(e)\|
<2\delta<\varepsilon,
\]
and
\[
\Lip(\beta)
\leq \delta\cdot \Lip(\varphi)+
\|K(e)\|\cdot \Lip(\psi)
\leq 2\delta+2\|K(e)\|
<4\delta<\varepsilon.
\]
Since \(\|T^{-1}\|\cdot \Lip(\beta)<\|T^{-1}\|\cdot \varepsilon<1\),
Lemma~\ref{lem:homeomorphism} implies that \(S\) is a homeomorphism.

Equip \(X/F\) with the quotient norm
\[
\|x+F\|_{X/F}:=\inf_{y\in F}\|x-y\|=\dist(x,F)\quad(x\in X).
\]
Let \(Q:X\to X/F\) be the canonical quotient map, given by \(Q(x)=x+F\).
Then \(Q\) is linear, continuous, and satisfies
\[
\|Q(x)\|_{X/F}=\dist(x,F)\leq\|x\|\quad(x\in X).
\]
Since \(\|\beta\|_{\infty}<\varepsilon<\varepsilon_{0}\) and
\(\Lip(\beta)<\varepsilon<\varepsilon_{0}\), Lemma~\ref{lem:orbit-condition}
applied to \(y=a_{0}\) yields \(x_{0}\in X\) such that
\(\|S^{m}(a_{0})-T^{m}(x_{0})\|\leq1\) for every \(m\in\Nzero\).
The identity \(S(a_{m})=a_{m+1}\) implies \(S^{m}(a_{0})=a_{m}\) for
\(m\in\{0,\ldots,N\}\), and hence
\begin{align*}
& \|T^{m}(x_{0})-a_{m}\|
=\|T^{m}(x_{0})-S^{m}(a_{0})\|\\
& \quad =\|S^{m}(a_{0})-T^{m}(x_{0})\|\leq1\quad(m\in\{0,\ldots,N\}).
\end{align*}
This, together with
\(\|a_{m}\|\leq 1+N\delta\), implies \(\|T^{m}(x_{0})\|\leq N\delta+2\) for
\(m\in\{0,\ldots,N\}\). Combining this with \eqref{eq:distance-KX-F}, we get
\[
\|Q(K(T^{m}(x_{0})))\|_{X/F}=
\dist(K(T^{m}(x_{0})),F)\leq\frac{\delta}{64}\quad(m\in\{0,\ldots,N\}).
\]
Let \(d_{m}=Q(T^{m}(x_{0}))-\lambda^{m}Q(x_{0})\)
for \(m\in\{0,\ldots,N\}\). Then \(d_{0}=Q(x_{0})-Q(x_{0})=0\).
Since \(T=\lambda I_{X}+K\) and \(Q\) is linear, we have
\(Q\circ T=\lambda Q+Q\circ K\). Thus, for every
\(m\in\{0,\ldots,N-1\}\),
\[
\begin{aligned}
d_{m+1}
&=Q(T^{m+1}(x_{0}))-\lambda^{m+1}Q(x_{0})\\
&=(Q\circ T)(T^{m}(x_{0}))-\lambda^{m+1}Q(x_{0})\\
&=\lambda Q(T^{m}(x_{0}))+Q(K(T^{m}(x_{0})))
-\lambda^{m+1}Q(x_{0})\\
&=\lambda(Q(T^{m}(x_{0}))-\lambda^{m}Q(x_{0}))
+Q(K(T^{m}(x_{0})))\\
&=\lambda d_{m}+Q(K(T^{m}(x_{0}))).
\end{aligned}
\]
Consequently,
\[
\|d_{m+1}\|_{X/F}\leq |\lambda|\cdot \|d_{m}\|_{X/F}
+\|Q(K(T^{m}(x_{0})))\|_{X/F}\leq\|d_{m}\|_{X/F}+\frac{\delta}{64}.
\]
Moreover, for every \(m\in\{1,\ldots,N\}\), we have
\[
\begin{aligned}
&\|Q(T^{m}(x_{0}))-\lambda^{m}Q(x_{0})\|_{X/F}
=\|d_{m}\|_{X/F}\leq\|d_{m-1}\|_{X/F}+\frac{\delta}{64}\\
&\quad\leq\cdots\leq\|d_{0}\|_{X/F}+\frac{m\delta}{64}=\frac{m\delta}{64}.
\end{aligned}
\]
Since \(\|Q\|\leq1\) and \(\|T^{m}(x_{0})-a_{m}\|\leq1\), we also have
\[
\|Q(T^{m}(x_{0}))-Q(a_{m})\|_{X/F}=
\|Q(T^{m}(x_{0})-a_{m})\|_{X/F}
\leq 1\quad(m\in\{0,\ldots,N\}).
\]
Combining these two inequalities and using the triangle inequality, we obtain
\[
\|\lambda^{m}Q(x_{0})-Q(a_{m})\|_{X/F}\leq1+\frac{m\delta}{64}\quad(m\in\{0,\ldots,N\}).
\]
This, together with \(Q(a_{m})=\lambda^{m}(1+m\delta)Q(e)\) and \(|\lambda|=1\), implies
\[
\|(1+m\delta)Q(e)-Q(x_{0})\|_{X/F}\leq1+\frac{m\delta}{64}\quad(m\in\{0,\ldots,N\}).
\]
In particular, taking \(m=0\) and \(m=N\) gives
\[
\|Q(e)-Q(x_{0})\|_{X/F}\leq 1
\quad\text{and}\quad
\|(1+N\delta)Q(e)-Q(x_{0})\|_{X/F}\leq1+\frac{N\delta}{64},
\]
implying
\[
N\delta\|Q(e)\|_{X/F}=\|((1+N\delta)Q(e)-Q(x_{0}))
-(Q(e)-Q(x_{0}))\|_{X/F}\leq2+\frac{N\delta}{64}.
\]
Since \(\|Q(e)\|_{X/F}=\dist(e,F)\geq\frac{1}{4}\), we get
\[
\frac{N\delta}{4}\leq2+\frac{N\delta}{64},
\quad\text{i.e.,}\quad\frac{15}{64}N\delta\leq2,
\]
which contradicts \(N\delta\geq16\).
Hence \(T\) is not strongly structurally stable.
\end{proof}

\begin{proposition}\label{prop:finite-spectrum}
Let \(X\) be a Banach space over \(\F\) and \(T\in\GL(X)\).
If there exists an isolated point \(\mu\) of \(\sigma(T)\) with
\(\mu\in\Tcircle\) such that the spectral subspace of \(T\) (if \(\F=\C\)) or
\(T_{\C}\) (if \(\F=\R\)) associated with \(\mu\) is finite-dimensional,
then \(T\) is not strongly structurally stable.
\end{proposition}

\begin{proof}
Define the Banach space \(Y\) over \(\F\), equipped with the absolute value norm, by
\[
Y=\begin{cases}
\R,&\F=\R\ \text{and}\ \mu\in\R,\\
\C,&\F=\C\ \text{or}\ \mu\in\C\setminus\R.
\end{cases}
\]

We first construct a bounded surjective \(\F\)-linear map
\(\Phi:X\to Y\) and a bounded \(\F\)-linear map
\(U:Y\to X\) such that
\[
\Phi\circ T=\mu\cdot\Phi \quad\text{and}\quad \Phi\circ U=I_{Y}.
\]

If \(\F=\C\), let \(P_{\mu}\) be the Riesz projection of \(T\)
associated with \(\mu\), and set \(E_{\mu}=P_{\mu}(X)\).
The space \(E_{\mu}\) is finite-dimensional,
\(T(E_{\mu})\subseteq E_{\mu}\), and \(\sigma(T|_{E_{\mu}})=\{\mu\}\)
\cite[Chapter~III, Theorem~6.17]{Kato95}.
Therefore there exists a nonzero functional \(\phi\in E_{\mu}^{*}\) such that
\[
\phi((T|_{E_{\mu}})(\xi))=\mu\phi(\xi)\quad(\xi\in E_{\mu}).
\]
Since \(P_{\mu}(X)=E_{\mu}\) and \(\phi\neq0\), the map
\(\Phi=\phi\circ P_{\mu}:X\to\C\) is nonzero, \(\C\)-linear, and bounded, with
\(\|\Phi\|\leq\|\phi\|\cdot \|P_{\mu}\|.\)
Since \(\Phi\) is a nonzero \(\C\)-linear functional, \(\Phi(X)=\C=Y\).
Choose \(u_{0}\in X\) with \(\Phi(u_{0})=1\) and define \(U:Y\to X\)
by \(U(y)=yu_{0}\).
The map \(U\) is \(\C\)-linear and satisfies
\[
\|U(y)\|=|y|\cdot \|u_{0}\|\quad(y\in Y),
\]
so \(U\) is bounded. Moreover,
\[
\Phi(U(y))=\Phi(yu_{0})=y\Phi(u_{0})=y\quad(y\in Y=\C),
\]
so \(\Phi\circ U=I_{Y}\).
Also, since \(P_{\mu}\circ T=T\circ P_{\mu}\), we have
\(\Phi\circ T=\phi\circ P_{\mu}\circ T
=\phi\circ(T|_{E_{\mu}})\circ P_{\mu}
=\mu\cdot(\phi\circ P_{\mu})=\mu\cdot\Phi.\)

If \(\F=\R\), let \(P_{\mu}\) be the Riesz projection of \(T_{\C}\) associated with
\(\mu\), and set \(E_{\mu}=P_{\mu}(X_{\C})\). The space \(E_{\mu}\) is
finite-dimensional, \(T_{\C}(E_{\mu})\subseteq E_{\mu}\), and
\(\sigma(T_{\C}|_{E_{\mu}})=\{\mu\}\)
\cite[Chapter~III, Theorem~6.17]{Kato95}.
Therefore there exists a nonzero functional \(\psi\in E_{\mu}^{*}\) such that
\[
\psi((T_{\C}|_{E_{\mu}})(\xi))=\mu\psi(\xi)\quad(\xi\in E_{\mu}).
\]
Since \(P_{\mu}(X_{\C})=E_{\mu}\) and \(\psi\neq0\), the map
\(\Psi=\psi\circ P_{\mu}:X_{\C}\to\C\) is nonzero, \(\C\)-linear, and bounded, with
\(\|\Psi\|\leq\|\psi\|\cdot \|P_{\mu}\|.\)
Also, since \(P_{\mu}\circ T_{\C}=T_{\C}\circ P_{\mu}\), we have
\(
\Psi\circ T_{\C}=\psi\circ P_{\mu}\circ T_{\C}
=\psi\circ(T_{\C}|_{E_{\mu}})\circ P_{\mu}
=\mu\cdot(\psi\circ P_{\mu})=\mu\cdot\Psi.
\)
The $\C$-linearity of \(\Psi\) gives
\[
\Psi(x_{1}+\mathbf{i}x_{2})=\Psi(x_{1})+\mathbf{i}\Psi(x_{2})
\quad(x_{1},x_{2}\in X).
\]
This, together with \(X_{\C}=X+\mathbf{i}X\) and \(\Psi\neq0\), implies
\(\Psi|_{X}\neq0\). We consider the following two cases:

If \(\mu\in\R\), then \(Y=\R\). Since \(\Psi|_{X}\neq0\), at least one of
\(\operatorname{Re}\circ(\Psi|_{X})\) and \(\operatorname{Im}\circ(\Psi|_{X})\) is nonzero.
Let \(\Phi:X\to\R\) be one such functional. Then \(\Phi\) is bounded and
\(\R\)-linear, with \(\|\Phi\|\leq\|\Psi\|\).
Since \(\Phi\) is a nonzero \(\R\)-linear functional, \(\Phi(X)=\R=Y\).
Choose \(u_{0}\in X\) with \(\Phi(u_{0})=1\) and define \(U:Y\to X\)
by \(U(y)=yu_{0}\).
The map \(U\) is \(\R\)-linear and satisfies
\[
\|U(y)\|=|y|\cdot \|u_{0}\|\quad(y\in Y),
\]
so \(U\) is bounded. Moreover,
\[
\Phi(U(y))=\Phi(yu_{0})=y\Phi(u_{0})=y\quad(y\in Y=\R),
\]
so \(\Phi\circ U=I_{Y}\).
Also, since \(T_{\C}|_{X}=T\) and \(\mu\in\R\), the identity
\(\Psi\circ T_{\C}=\mu\cdot\Psi\) gives, for every \(x\in X\),
\[
\operatorname{Re}(\Psi(T(x)))
=\operatorname{Re}(\mu\Psi(x))
=\mu\operatorname{Re}(\Psi(x)),
\]
and
\[
\operatorname{Im}(\Psi(T(x)))
=\operatorname{Im}(\mu\Psi(x))=\mu\operatorname{Im}(\Psi(x)).
\]
Thus \(\Phi\circ T=\mu\cdot\Phi\).

If \(\mu\in\C\setminus\R\), then \(Y=\C\). Set
\(\Phi=\Psi|_{X}:X\to\C\). The map \(\Phi\) is nonzero, \(\R\)-linear,
and bounded, with \(\|\Phi\|\leq\|\Psi\|\). Also, since \(T_{\C}|_{X}=T\)
and \(\Psi\circ T_{\C}=\mu\cdot\Psi\), for every \(x\in X\), we have
\(\Phi(T(x))=\Psi(T(x))=\Psi(T_{\C}(x))=\mu\Psi(x)=\mu\Phi(x).\)
Since \(T\in\GL(X)\), we have \(T(X)=X\). This, together with
\(\Phi\circ T=\mu\cdot\Phi\) and \(\Phi\neq0\), gives
\(\mu\Phi(X)=\Phi(T(X))=\Phi(X)\neq \{0\}.\)
Choose \(y_{0}\in\Phi(X)\setminus\{0\}\). Then \(\mu y_{0}\in\Phi(X)\).
Since \(y_{0}\neq0\) and \(\operatorname{Im}(\mu)\neq0\), for all \(a,b\in\R\),
\[
ay_{0}+b\mu y_{0}=0
\ \Longrightarrow\ a+b\mu=0
\ \Longrightarrow\ a=b=0.
\]
Thus \(y_{0}\) and \(\mu y_{0}\) are linearly independent over \(\R\), and hence
\(\C=\spanop\{y_{0},\) \(\mu y_{0}\}\subseteq\Phi(X)\subseteq\C.\)
Therefore \(\Phi(X)=\C=Y\), so \(\Phi\) is surjective.
Choose \(u_{0},v_{0}\in X\) with \(\Phi(u_{0})=1\) and \(\Phi(v_{0})=\mathbf{i}\), and define
\(U(a+\mathbf{i}b)=au_{0}+bv_{0}\) for all \(a,b\in\R\). Then \(U\) is \(\R\)-linear and
\[
\Phi(U(a+\mathbf{i}b))
=\Phi(au_{0}+bv_{0})
=a\Phi(u_{0})+b\Phi(v_{0})=a+\mathbf{i}b\quad(a,b\in\R),
\]
so \(\Phi\circ U=I_{Y}\). Moreover,
\(\|U(a+\mathbf{i}b)\|\leq(\|u_{0}\|+\|v_{0}\|)|a+\mathbf{i}b|\), so \(U\) is bounded.

In all cases for $Y$, the maps \(\Phi\) and \(U\) satisfy
\(\Phi\circ T=\mu\cdot\Phi\) and \(\Phi\circ U=I_{Y}\).

Suppose, on the contrary, that \(T\) is strongly structurally stable. Let
\(\varepsilon_{0}>0\) be given by Lemma \ref{lem:orbit-condition} with
\(\gamma=1\), and set
\[
\varepsilon:=\frac{1}{2}\min\left\{1,\|T^{-1}\|^{-1},\varepsilon_{0}\right\},\quad
\delta:=\frac{\varepsilon}{4\|U\|},
\quad\text{and}\quad N:=\left\lfloor\frac{2\|\Phi\|}{\delta}\right\rfloor+1.
\]
Then \(N\delta>2\|\Phi\|\).
Define
\[
a_{0}=U(\|\Phi\|) \quad\text{and}\quad a_{m+1}=T(a_{m})+\delta U(\mu^{m+1})
\quad(m\in\{0,\ldots,N-1\}).
\]
The identities \(\Phi\circ T=\mu\cdot\Phi\) and \(\Phi\circ U=I_{Y}\) give
\(\Phi(a_{0})=\Phi(U(\|\Phi\|))=\|\Phi\|\) and
\begin{align*}
&\Phi(a_{m+1})
=\Phi(T(a_{m}))+\delta\Phi(U(\mu^{m+1}))\\
&\quad =\mu\Phi(a_{m})+\delta\mu^{m+1}\quad (m\in\{0,\ldots,N-1\}).
\end{align*}
Moreover, for every \(m\in\{1,\ldots,N\}\), we have
\[
\Phi(a_{m})
=\mu\Phi(a_{m-1})+\delta\mu^{m}
=\cdots=\mu^{m}\Phi(a_{0})+m\delta\mu^{m}
=\mu^{m}(\|\Phi\|+m\delta).
\]

The assignment
\[
a_{m}\mapsto\mu^{m+1}\quad(m\in\{0,\ldots,N\})
\]
is well defined because \(|\Phi(a_{m})|=\|\Phi\|+m\delta\) shows that
\(a_{0},\ldots,a_{N}\) are pairwise distinct. Lemma
\ref{lem:radial}, applied with \(\nu=\mu\) and \(R=\|\Phi\|\), gives
\[
|\mu^{m+1}-\mu^{n+1}|
\leq \frac{1}{\|\Phi\|}|\Phi(a_{m})-\Phi(a_{n})|
\leq\|a_{m}-a_{n}\|
\]
for \(m,n\in\{0,\ldots,N\}\). Hence the assignment has
Lipschitz constant at most \(1\).

Since \(|\mu|=1\), for every \(m\in\{0,\ldots,N\}\), we have
\(|\mu^{m+1}|=1.\)
Lemma \ref{lem:interpolation}, applied with \(R=1\),
therefore gives a bounded Lipschitz map \(G:X\to Y\) such that
\[
G(a_{m})=\mu^{m+1}\quad(m\in\{0,\ldots,N\}),
\]
and
\[
\|G\|_{\infty}\leq1,\quad \Lip(G)\leq2 .
\]

Define
\(\beta(x)=\delta U(G(x)) .\)
The interpolation identities for \(G\), together with
\(a_{m+1}=T(a_{m})+\delta U(\mu^{m+1})\), give
\[
\beta(a_{m})
=\delta U(G(a_{m}))
=\delta U(\mu^{m+1})
=a_{m+1}-T(a_{m})\quad(m\in\{0,\ldots,N-1\}).
\]
Thus the map \(S=T+\beta\) satisfies
\[
S(a_{m})=a_{m+1}\quad(m\in\{0,\ldots,N-1\}).
\]
Moreover, the choice of \(\delta\) and the bounds on \(G\) give
\[
\|\beta\|_{\infty}\leq\delta\cdot\|U\|\cdot\|G\|_{\infty}
\leq\frac{\varepsilon}{4}<\varepsilon,
\]
and
\[
\Lip(\beta)\leq\delta\cdot\|U\|\cdot\Lip(G)
\leq\frac{\varepsilon}{2}<\varepsilon.
\]
Since \(\|T^{-1}\|\cdot\Lip(\beta)<\|T^{-1}\|\cdot\varepsilon<1\),
Lemma \ref{lem:homeomorphism} implies that \(S\) is a homeomorphism.

Since \(\|\beta\|_{\infty}<\varepsilon<\varepsilon_{0}\) and
\(\Lip(\beta)<\varepsilon<\varepsilon_{0}\), Lemma \ref{lem:orbit-condition}
applied to \(y=a_{0}\) yields \(x_{0}\in X\) such that
\(\|S^{m}(a_{0})-T^{m}(x_{0})\|\leq1\) for every \(m\in\Nzero\).
The identities \(S(a_{m})=a_{m+1}\) imply \(S^{m}(a_{0})=a_{m}\) for
\(m\in\{0,\ldots,N\}\), and hence
\begin{align*}
& \|T^{m}(x_{0})-a_{m}\|
=\|T^{m}(x_{0})-S^{m}(a_{0})\|\\
&\quad =\|S^{m}(a_{0})-T^{m}(x_{0})\|\leq 1
\quad(m\in\{0,\ldots,N\}).
\end{align*}
Iterating \(\Phi\circ T=\mu\cdot\Phi\) gives
\(\Phi(T^{m}(x_{0}))=\mu^{m}\Phi(x_{0})\) for every \(m\in\Nzero\).
Since \(\Phi\) is bounded and \(\|T^{m}(x_{0})-a_{m}\|\leq1\), we also have
\[
\begin{aligned}
& |\Phi(T^{m}(x_{0}))-\Phi(a_{m})|
=|\Phi(T^{m}(x_{0})-a_{m})|\\
&\quad \leq\|\Phi\|\cdot\|T^{m}(x_{0})-a_{m}\|
\leq\|\Phi\|\quad(m\in\{0,\ldots,N\}).
\end{aligned}
\]
Since \(\Phi(T^{m}(x_{0}))=\mu^{m}\Phi(x_{0})\),
\(\Phi(a_{m})=\mu^{m}(\|\Phi\|+m\delta)\), and \(|\mu|=1\),
the preceding estimate yields
\[
\begin{aligned}
& \bigl|\|\Phi\|+m\delta-\Phi(x_{0})\bigr|
=\bigl|\mu^{m}(\Phi(x_{0})-\|\Phi\|-m\delta)\bigr|\\
&\quad =|\Phi(T^{m}(x_{0}))-\Phi(a_{m})|
\leq\|\Phi\|\quad(m\in\{0,\ldots,N\}).
\end{aligned}
\]
In particular, taking \(m=0\) and \(m=N\) gives
\[
\bigl|\|\Phi\|-\Phi(x_{0})\bigr|\leq\|\Phi\|
\quad\text{and}\quad
\bigl|\|\Phi\|+N\delta-\Phi(x_{0})\bigr|\leq\|\Phi\|,
\]
implying
\(N\delta\leq2\|\Phi\|,\)
which contradicts \(N\delta>2\|\Phi\|\). Hence \(T\) is not strongly
structurally stable.
\end{proof}

We point out that the proofs of Propositions~\ref{prop:essential-scalar} and
\ref{prop:finite-spectrum} use Lemmas~\ref{lem:radial} and
\ref{lem:interpolation} to construct small bounded Lipschitz perturbations with
prescribed finite orbit segments. After applying a bounded linear map, the norm
along each prescribed segment grows linearly, leading to a contradiction with
the orbit comparison in Lemma~\ref{lem:orbit-condition}.
The difference lies in the treatment of the unperturbed orbits.
In Proposition~\ref{prop:essential-scalar}, the identity
\(Q\circ T=\lambda\cdot Q+Q\circ K\) contains the error term \(Q\circ K\),
which is controlled on the relevant bounded set by compactness.
In Proposition~\ref{prop:finite-spectrum}, the finite-dimensional spectral subspace
instead yields the exact identity \(\Phi\circ T=\mu\cdot\Phi\).

\begin{proof}[Proof of Theorem \ref{thm:main}]
Suppose, on the contrary, that \(T\) is not hyperbolic,
i.e., \(\sigma(T)\cap\Tcircle\neq\varnothing\). Then,
fix \(\mu\in\sigma(T)\cap\Tcircle\).

We first claim that \(X\) is infinite-dimensional.
Indeed, if \(X\) is finite-dimensional, then \(T\) (if \(\F=\C\)) or
\(T_{\C}\) (if \(\F=\R\)) acts on a finite-dimensional complex Banach space.
Hence \(\mu\) is an isolated point of \(\sigma(T)\), and its associated spectral
subspace for the corresponding complex operator is finite-dimensional.
Proposition~\ref{prop:finite-spectrum}, applied to \(\mu\), then contradicts the strong structural stability of
\(T\). Thus \(X\) is infinite-dimensional. This, together with
Proposition~\ref{prop:essential-scalar}, implies \(|\lambda|\neq 1\).

Since \(|\mu|=1\) and \(|\lambda|\neq1\), we have \(\mu\neq\lambda\).
To apply Proposition~\ref{prop:finite-spectrum}, we now show that \(\mu\) is an
isolated point of \(\sigma(T)\) and that the spectral subspace of
\(T\) (if \(\F=\C\)) or \(T_{\C}\) (if \(\F=\R\)) associated with \(\mu\) is
finite-dimensional. We verify these two properties separately in the
following two cases:

\textbf{Case~1.} \(\F=\C\). The identity \(\sigma(T)=\lambda+\sigma(K)\) gives
\(\mu-\lambda\in\sigma(K)\setminus\{0\}.\)
By the Riesz--Schauder theorem~\cite[Chapter~III, Theorem~6.26]{Kato95},
applied to the compact operator \(K\), the point \(\mu-\lambda\) is isolated
in \(\sigma(K)\) and its associated spectral subspace is finite-dimensional.
Since \(T=\lambda I_{X}+K\), the point \(\mu\) is isolated in \(\sigma(T)\),
and the Riesz projection of \(T\) associated with \(\mu\) equals that of
\(K\) associated with \(\mu-\lambda\). Hence the spectral subspace of
\(T\) associated with \(\mu\) is finite-dimensional.

\textbf{Case~2.} \(\F=\R\). We have \(T_{\C}=\lambda I_{X_{\C}}+K_{\C}\).
Let \(B\subseteq X_{\C}\) be an arbitrary bounded set.
The real and imaginary parts of its elements form bounded subsets of \(X\),
whose images under \(K\) are relatively compact. The identity
\(K_{\C}(x+\mathbf{i}y)=K(x)+\mathbf{i}K(y)\) therefore shows that
\(K_{\C}(B)\) is relatively compact. Thus \(K_{\C}\) is compact.
Since \(\sigma(T)=\sigma(T_{\C})=\lambda+\sigma(K_{\C})\), we have
\(\mu-\lambda\in\sigma(K_{\C})\setminus\{0\}.\)
By the Riesz--Schauder theorem~\cite[Chapter~III, Theorem~6.26]{Kato95},
applied to the compact operator \(K_{\C}\), the point \(\mu-\lambda\)
is isolated in \(\sigma(K_{\C})\)
and its associated spectral subspace is finite-dimensional.
Since \(T_{\C}=\lambda I_{X_{\C}}+K_{\C}\), the point \(\mu\) is isolated
in \(\sigma(T_{\C})\), and the Riesz projection of \(T_{\C}\) associated
with \(\mu\) equals that of \(K_{\C}\) associated with \(\mu-\lambda\).
Hence the spectral subspace of \(T_{\C}\) associated with \(\mu\) is
finite-dimensional.

In both cases, since \(\mu\in\Tcircle\), Proposition~\ref{prop:finite-spectrum}
implies that \(T\) is not strongly structurally
stable, a contradiction.

Therefore \(\sigma(T)\cap\Tcircle=\varnothing\), so \(T\) is hyperbolic.
\end{proof}

\begin{remark}\label{rem:negative}
By Theorem~\ref{thm:argyros-haydon} and Theorem~\ref{thm:ah-complexification},
respectively, the real space \(\AHX\) and
the complex space \((\AHX)_{\C}\) are infinite-dimensional separable Banach
spaces with the scalar-plus-compact property. Thus every bounded linear
operator on either space has the form \(\lambda I+K\), where \(\lambda\)
belongs to the corresponding scalar field and \(K\) is compact.
Theorem~\ref{thm:main} therefore implies that every strongly structurally
stable operator on either space is hyperbolic.
Thus these infinite-dimensional separable real and complex Banach spaces
give a negative answer to the part of Question~\ref{q:bm} concerning strongly
structurally stable operators.
\end{remark}

\section*{Acknowledgements}
We sincerely thank the reviewers for their insightful comments and constructive suggestions, which have greatly helped improve this manuscript.

\end{document}